\documentclass[12pt, reqno]{amsart}
\usepackage{graphicx}

\newtheorem{theorem}{Theorem}[section]
\newtheorem{proposition}[theorem]{Proposition}
\newtheorem{lemma}[theorem]{Lemma}

\theoremstyle{definition} 
\newtheorem{definition}[theorem]{Definition}

\newtheorem{remark}[theorem]{Remark}

\newcommand{\C}{\mathbb C} \newcommand{\N}{\mathbb N}
\newcommand{\R}{\mathbb R} 
\newcommand{\Z}{\mathbb Z}

\DeclareMathOperator{\dist}{{\rm dist}}
\DeclareMathOperator{\sys}{{\rm sys}}

\DeclareMathOperator{\pisys}{\sys\!\pi}

\DeclareMathOperator{\area}{{\rm area}}
\DeclareMathOperator{\length}{{\rm length}}

\DeclareMathOperator{\Aut}{{\rm Aut}}

\DeclareMathOperator{\gmetric}{{\mathcal G}}
\DeclareMathOperator{\Bolza}{{\mathcal B}}

\DeclareMathOperator{\CAT}{{\rm CAT}}

\def\SR{{\rm SR}}

\def\go {{\gmetric_{\mathcal O}}} \def\genus {{\it g}}

\def\ie {{\it i.e.\ }} 
 
\def\cf {\hbox{\it cf.\ }}
\def\rat{{\angle\!^\mid \! \left(\frac{x}{2}, \frac{\pi}{8}\right)\ }}

\long\def\forget#1\forgotten{} %

\newcommand{\rp}{\mathbb R\mathbb P} 

\numberwithin{equation}{section} 
\numberwithin{figure}{section}

\begin{document}

\author[M.~Katz]{Mikhail G. Katz$^{*}$} \address{Department of
Mathematics, Bar Ilan University, Ramat Gan 52900 Israel}
\email{katzmik@math.biu.ac.il} \thanks{$^{*}$Supported by the Israel
Science Foundation (grants no.\ 620/00-10.0 and 84/03)}

\author[S.~Sabourau]{St\'ephane Sabourau} \address{Laboratoire de
Math\'ematiques et Physique Th\'eorique, Universit\'e de Tours, Parc
de Grandmont, 37400 Tours, France}

\address{Department of Mathematics, University of Pennsylvania,
209 South 33rd Street, Philadelphia, PA 19104-6395, USA}
\email{sabourau@lmpt.math.univ-tours.fr}


\title[An optimal inequality for CAT(0) metrics in genus two]{An
optimal systolic inequality for CAT(0) metrics in genus two}


\begin{abstract}
We prove an optimal systolic inequality for CAT(0) metrics on a
genus~2 surface.  We use a Voronoi cell technique, introduced by
C.~Bavard in the hyperbolic context.  The equality is saturated by a
flat singular metric in the conformal class defined by the smooth
completion of the curve $y^2=x^5-x$.  Thus, among all CAT(0) metrics,
the one with the best systolic ratio is composed of six flat regular
octagons centered at the Weierstrass points of the Bolza surface.
\end{abstract}

\maketitle

\tableofcontents

\section{Hyperelliptic surfaces of nonpositive curvature}
\label{Jone1}

Over half a century ago, a student of C. Loewner's named P. Pu
presented, in the pages of the Pacific Journal of Mathematics
\cite{Pu}, the first two optimal systolic inequalities, which came to
be known as the Loewner inequality for the torus, and Pu's
inequality~\eqref{J62} for the real projective plane.

The recent months have seen the discovery of a number of new systolic
inequalities \cite{Am, BK1, Sa1, BK2, IK, BCIK1, BCIK2, KL, Ka4, KS1,
KRS}, as well as near-optimal asymptotic bounds \cite{Ka3, KS2, Sa2,
KSV, Sa3, RS}.  A number of questions posed in \cite{CK} have thus
been answered.  A general framework for systolic geometry in a
topological context is proposed in~\cite{KR, KR2}.

The homotopy 1-systole, denoted $\pisys_1(X)$, of a compact metric
space $X$ is the least length of a noncontractible loop of $X$.

Given a metric~$\gmetric$ on a surface, let~$\SR(\gmetric)$ denote its
systolic ratio
\[
\SR(\gmetric) =\frac{\pisys_1(\gmetric)^2}{\area(\gmetric)}.
\]
Given a compact Riemann surface~$\Sigma$, its {\em optimal systolic
ratio\/} is defined as~$\SR(\Sigma)= \sup_{\gmetric} \SR(\gmetric)$,
where the supremum is over all metrics in the conformal type
of~$\Sigma$.  Finally, given a smooth compact surface~$M$, its optimal
systolic ratio is defined by setting~$\SR(M)= \sup_\Sigma
\SR(\Sigma)$, where the supremum is over all conformal
structures~$\Sigma$ on~$M$.  The latter ratio is known for the Klein
bottle, \cf \eqref{Jcb}, in addition to the torus and real projective
plane already mentioned.

In the class of all metrics without any curvature restrictions, no
singular flat metric on a surface of genus~2 can give the optimal
systolic ratio in this genus \cite{Sa1}.  The best available upper
bound for the systolic ratio of an arbitrary genus two surface
is~$\gamma_2 \simeq 1.1547$, \cf \cite{KS1}.

The precise value of~$\SR$ for the genus~2 surface has so far eluded
researchers, \cf~\cite{Ca, Bry}.  We propose an answer in the
framework of negatively curved, or more generally,~$\CAT(0)$ metrics.

The term ``$\CAT(0)$ space'' evokes an extension of the notion of a
manifold of nonpositive curvature, to encompass singular spaces.  We
will use the term to refer to surfaces with metrics with only mild
quotient singularities, defined below.  Here the condition of
nonpositive curvature translates into a lower bound of~$2\pi$ for the
total angle at the singularity.  We need such an extension so as to
encompass the metric which saturates our optimal
inequality~\eqref{J11}.

A mild quotient singularity is defined as follows.  Consider a smooth
metric on~$\R^2$.  Let~$q\geq 1$ be an integer.  Consider the~$q$-fold
cover~$X_q$ of~$\R^2\setminus \{0\}$ with the induced metric.  We
compactify~$X_q$ in the neighborhood of the origin to obtain a
complete metric space~$X_q^c = X_q \cup \{0\}$.  

\begin{definition}
\label{911}
Suppose~$X_q^c$ admits an isometric action of~$\Z_p$ fixing the
origin.  Then we can form the orbit space~$Y_{p,q} = X_q^c/\Z_p$.  The
space~$Y_{p,q}$ is then called mildly singular at the origin.  
\end{definition}

The total angle at the singularity is then $\tfrac {2\pi q}{p}$, and
the CAT(0) condition is $\tfrac q p \geq 1$.

\begin{remark}
\label{tro}
Alternatively, a point is singular of total angle $2\pi(1+\beta)$ if
the metric is of the form $e^{h(z)} |z|^{2\beta} |dz|^2$ in its
neighborhood, where $|dz|^2=dx^2 + dy^2$, \cf
M.~Troyanov~\cite[p.~915]{Tro}.
\end{remark}

\begin{theorem} 
\label{Jtheo:nonpositive}
Every~$\CAT(0)$ metric~$\gmetric$ on a surface~$\Sigma_2$ of genus
two, satisfies the following optimal inequality:
\begin{equation}
\label{J11}
\SR(\Sigma_2,\gmetric) \leq \tfrac{1}{3} \cot\tfrac{\pi}{8} =
\tfrac{1}{3}\left( \sqrt{2} + 1 \right) = 0.8047\ldots
\end{equation}
The inequality is saturated by a singular flat metric, with~$16$
conical singularities, in the conformal class of the Bolza
surface.
\end{theorem}

The Bolza surface is described in Section~\ref{J22}.  The optimal
metric is described in more detail in Section~\ref{Jgenus2}.
Theorem~\ref {Jtheo:nonpositive} is proved in Section~\ref{Jproof12}
based on the octahedral triangulation of~$S^2$.

\begin{remark}
A similar optimal inequality can be proved for hyperelliptic surfaces
of genus~5 based on the icosahedral triangulation, \cf~\cite{Bav1}.
\end{remark}

\section{Distinguishing $16$ points on the Bolza surface}
\label{J22}

The Bolza surface~$\Bolza$ is the smooth completion of the affine
algebraic curve
\begin{equation}
\label{J71}
y^2= x^5-x.
\end{equation}
It is the unique Riemann surface of genus two with a group of
holomorphic automorphisms of order 48.  (A way of passing from an
affine hyperelliptic surface to its smooth completion is described in
\cite[p.~60-61]{Mi}.)

\begin{definition}
A conformal involution~$J$ of a compact Riemann surface~$\Sigma$ of
genus $\genus$ is called {\em hyperelliptic\/} if~$J$ has
precisely~$2\genus + 2$ fixed points.  The fixed points of~$J$ are
called the Weierstrass points of~$\Sigma$.
\end{definition}

The quotient Riemann surface~$\Sigma/J$ is then necessarily the
Riemann sphere, denoted henceforth~$S^2$.  Let~$Q: \Sigma \to S^2$ be
the conformal ramified double cover, with~$2\genus+2$ branch points.
Thus,~$J$ acts on~$\Sigma$ by sheet interchange.  Recall that every
surface of genus~2 is hyperelliptic, \ie admits a hyperelliptic
involution \cite[Proposition III.7.2]{FK}.

We make note of $16$ special points on~$\Bolza$.  We refer to a point
as being special if it is a fixed point of an order 3 automorphism
of~$\Bolza$.

Consider the regular octahedral triangulation of~$S^2=\C\cup \infty$.
Its set of vertices is conformal to the set of roots of the polynomial
$x^5-x$ of formula~\eqref{J71}, together with the unique point at
infinity.  Thus the six points in question can be thought of as the
ramification points of the ramified conformal double cover~$Q:{\Bolza}
\to S^2$, while the $16$ special points of~$\Bolza$ project to the 8
vertices of the cubical subdivision dual to the octahedral
triangulation.

In other words, the~$x$-coordinates of the ramification points are
\[
\left\{ 0, \infty, 1, -1, i, -i \right\},
\]
which stereographically correspond to the vertices of a regular
inscribed octahedron.  The conformal type therefore admits the
symmetries of the cube.  If one includes both the hyperelliptic
involution, and the real (antiholomorphic) involution of~$\Bolza$
corresponding to the complex conjugation~$(x,y)\to (\bar{x}, \bar{y})$
of~$\C ^2$, one obtains the full symmetry group~$\Aut(\Bolza)$, of
order
\begin{equation}
\label{J51}
|\Aut(\Bolza)| = 96,
\end{equation}
see \cite[p.~404]{KuN} for more details.

\begin{lemma}
\label{J21}
The hyperbolic metric of~$\Bolza$ admits~$12$~systolic loops.
The~$12$ loops are in one to one correspondence with the edges of the
octahedral decomposition of~$S^{2}$.  The correspondence is given by
taking the inverse image under~$Q$ of an edge.  The~$12$~systolic
loops cut the surface up into~$16$~hyperbolic triangles.  The centers
of the triangles are the~$16$~special points.
\end{lemma}

See \cite[\S~5]{Sc1} for further details.  The Bolza surface is
extremal for two distinct problems:
\begin{itemize}
\item
systole of hyperbolic surfaces \cite{Bav2}, \cite[Theorem 5.2]{Sc1};
\item
conformal systole of Riemann surfaces \cite{BS}.
\end{itemize}
The square of the conformal systole of a Riemann surface is also known
as its Seshadri constant \cite{Kon}.  The Bolza surface is also
conjectured to be extremal for the first eigenvalue of the Laplacian.
Such extremality has been verified numerically \cite{JP}.  The above
evidence suggests that the systolically extremal surface may lie in
the conformal class of~$\Bolza$, as well.  Meanwhile, we have the
following result, proved in Section~\ref{Jfive}.

\begin{theorem}
\label{J13}
The Bolza surface~$\Bolza$ satisfies~$\SR(\Bolza) \leq \frac{\pi}{3}$.
\end{theorem}
Note that Theorems~\ref{J13} and \ref{Jtheo:nonpositive} imply that
$\SR(\Bolza) \in [0.8, 1.05]$.

\section{A flat singular metric in genus~two}
\label{Jgenus2}

The optimal systolic ratio of a genus~2 surface~$(\Sigma_2, \gmetric)$
is unknown, but it satisfies the Loewner inequality, \cf \cite{KS1}.
Here we discuss a {\em lower\/} bound for the optimal systolic ratio
in genus~2, briefly described in~\cite{CK}.

The example of M. Berger (see \cite [Example~5.6.B$'$]{Gr1}) in
genus~2 is a singular flat metric with conical singularities.  It has
systolic ratio~$\SR=0.6666$.  This ratio was improved by F. Jenni
\cite{Je}, who identified the hyperbolic genus~2 surface with the
optimal systolic ratio among all hyperbolic genus~2 surfaces (see also
C.~Bavard~\cite{Bav2} and P.~Schmutz \cite[Theorem~5.2]{Sc1}).  The
surface in question is a triangle surface~(2,3,8).  Its conformal
class is that of the Bolza surface, \cf Section~\ref{J22}.  It admits
a regular hyperbolic octagon as a fundamental domain, and has
12~systolic loops of length~$2x$, where~$x=\cosh^{-1} (1+\sqrt{2})$.
It has
\[
\pisys_1= 2\log\left( 1+\sqrt{2}+ \sqrt{2+2\sqrt{2}} \right),
\]
area~${4\pi}$, and systolic ratio~$\SR\simeq 0.7437$.  This ratio can
be improved as follows (see the table of Figure~\ref{Jfig2}).

\begin{figure}
\renewcommand{\arraystretch}{1.8} 
\[
\begin{tabular}[t]{|
@{\hspace{3pt}}p{1.5in}|| @{\hspace{3pt}}p{.8in}|
@{\hspace{3pt}}p{.8in}| @{\hspace{3pt}}p{.8in}| } \hline & Berger &
Jenni & metric $\go$ on Bolza \\ \hline\hline $\displaystyle
\SR(\gmetric) =\frac{\pisys_1(\gmetric)^2}{\area(\gmetric)}$ & 0.6666
& 0.7437 & 0.8047 \\ \hline
\end{tabular}
\]
\renewcommand{\arraystretch}{1}
\caption{\textsf{Three systolic ratios in genus~2}}
\label{Jfig2}
\end{figure}

\begin{proposition}
The conformal class of the Bolza surface~$\Bolza$ admits a metric,
denoted~$\go$, with the following properties:
\begin{enumerate}
\item
the metric is singular flat, with conical singularities precisely at
the $16$ special points of Section~$\ref{J22}$;
\item
each singularity is of total angle~$\frac{9\pi}{4}$, so that the
metric~$\go$ is $\CAT(0)$;
\item
the metric is glued from six flat regular octagons, centered on the
Weierstrass points, while the 1-skeleton projects under~$Q:{\Bolza}
\to S^2$ to that of the dual cube in~$S^2$;
\item
the systolic ratio equals~$\SR(\go)= \tfrac{1}{3}\left( \sqrt{2} + 1
\right) > \frac{4}{5}$.
\end{enumerate}
\end{proposition}

\begin{figure}
\setlength\unitlength{1pt}
\noindent \begin{picture}
(400,150)(0,0)
\put(120,15){\includegraphics[angle=0,height=120pt]{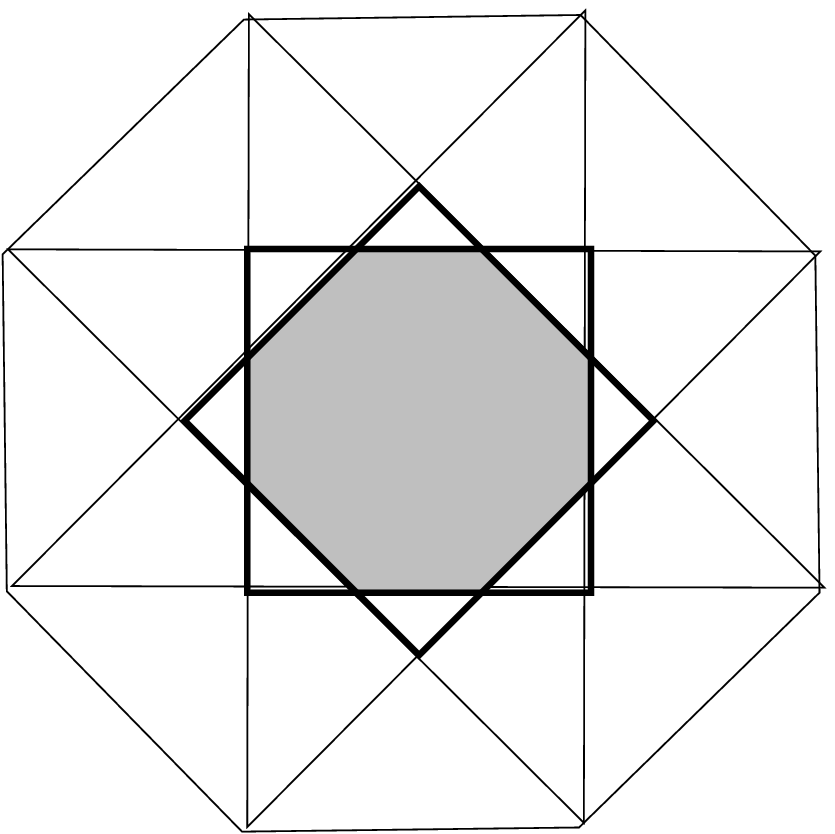}}
\end{picture}
\caption{\textsf{Flat regular octagon obtained as the union of $16$
triangles $\rat$.  The shaded interior octagon represents the region
with 4 geodesic loops through every point}}
\label{Jfigure1}
\end{figure}

\begin{proof} 
The octahedral triangulation of the sphere, \cf Section~\ref{J22},
lifts to a triangulation of~$\Bolza$ consisting of $16$ triangles,
which we think of as being ``equilateral''.  Here 8 equilateral
triangles are connected cyclically around each of the 6 Weierstrass
vertices of the triangulation of~$\Bolza$.

We further subdivide each equilateral triangle into 3 isosceles
triangles, with a common vertex at the center of the equilateral
triangle.  We equip each of the 48 isosceles triangles with a flat
metric with obtuse angle $\tfrac{3\pi}{4}$.

Each of the 6 Weierstrass vertices of the original triangulation is a
smooth point, since the total angle is ~$8\left(\tfrac{\pi}{4} \right)
=2\pi$.  Each equilateral triangle possesses a singularity at the
center with total angle~$\tfrac{9\pi}{4}=2\pi\left( 1+\tfrac{1}{8}
\right)> 2\pi$.  Alternatively, we can apply the Gauss-Bonnet formula
$\sum_{\sigma} \alpha (\sigma) = 2\genus-2$, in genus 2 with $16$
isometric singularities.  Here the sum is over all
singularities~$\sigma$ of a singular flat metric on a surface of
genus~$\genus$, where the cone angle at singularity~$\sigma$
is~$2\pi(1+\alpha(\sigma))$.  Since the metric~$\go$ is smooth at a
Weierstrass point of~$\Bolza$, the metric has
only $16$ singularities, precisely at the special points of
Section~\ref{J22}, proving items 1 and 2 of the proposition.

Let $x$ denote the sidelength of the equilateral triangle.  The
barycentric subdivision of each equilateral triangle consists of 6
copies of a flat right angle triangle, denoted $\rat$, with
side~$\frac{x}{2}$ and adjacent angle~$\frac{\pi}{8}$.  We thus obtain
a decomposition of the metric~$\go$ into 96 copies of the
triangle~$\rat$, which can be thought of as a fundamental domain for
the action of~$\Aut(\Bolza)$, \cf \eqref{J51}.


We have $\pisys_1(\go)= 2x$ by Lemma~\ref{J31}, proving item 4 of the
proposition.  To prove item 3, note that the union of the $16$
triangles~$\rat$ with a common Weierstrass vertex is a flat regular
octagon.  The latter is represented in Figure~\ref{Jfigure1}, together
with the systolic loops passing through~it.
\end{proof}

\begin{lemma}
\label{J31}
The systole of the singular flat $\CAT(0)$ metric on the Bolza
surface equals twice the distance between a pair
of adjacent Weierstrass points.
\end{lemma}

\begin{proof}
Consider the smooth closed geodesic $\gamma \subset \Bolza$ which is
the inverse image under the map $Q: \Bolza \to S^2$ of an edge of the
octahedron, \cf Lemma~\ref{J21}.  Let~$x$ be the distance between a
pair of opposite sides of the regular flat octagon in $\Bolza$, or
equivalently, the distance between a pair of adjacent Weierstrass
points.  Thus, $\length(\gamma)=2x$.  Consider a loop~$\alpha \subset
\Bolza$ whose length satisfies
\begin{equation}
\label{941}
\length(\alpha) < 2x.
\end{equation}
We will prove that there are two possibilities for $\alpha$: it is
either contractible, or freely homotopic to one of the 12 geodesics
$\gamma$ of the type described above.  On the other hand, $\gamma$ is
necessarily length minimizing in its free homotopy class, by the
$\CAT(0)$ property of the metric \cite[Theorem~6.8]{BH}.  This will
rule out the second possibility, and prove the lemma.

Denote by~$\Bolza^{(1)}$ the graph on~$\Bolza$ given by the inverse
image under~$Q$ of the 1-skeleton of the cubical subdivision of~$S^2$.
The graph~$\Bolza^{(1)}$ partitions the surface into six regular
octagons, denoted $\Omega_k$:
\begin{equation}
\label{942}
\Bolza= \cup_{k=1}^6 \Omega_k.
\end{equation}
We will deform $\alpha$ to a loop~$\beta \subset \Bolza^{(1)}$, as
follows.  The partition \eqref{942} induces a partition of the
loop~$\alpha$ into arcs~$\alpha_i$, each lying in its respective
octagon $\Omega_{k_i}$.  We deform each~$\alpha_i$, without increasing
length, to the straight line segment $[p_i, q_i]\subset \Omega_{k_i}$.
The boundary points $p_i, q_i$ of~$\alpha_i$ split the boundary closed
curve $\partial \Omega_{k_i} \subset \Bolza^{(1)}$ into a pair of
paths.  Let~$\beta_i \subset \partial \Omega_{k_i}$ be the shorter of
the two paths.  Denote by~$y$ the distance between adjacent vertices
of the octagon.  Then clearly
\begin{equation}
\label{J32}
\length (\beta_i) \leq \frac{4y}{x} \length(\alpha_i).
\end{equation}
We first deform the loop $\alpha$ into the graph~$\Bolza^{(1)}$.  The
deformation fixes the intersection points~$\alpha \cap \Bolza^{(1)}$.
Inside $\Omega_{k_i}$, we deform the arc~$\alpha_i$ to the
path~$\beta_i$.
The length of the resulting loop is at most
\[
\frac{4y}{x} \length(\alpha) < 8 y
\]
by \eqref{941} and \eqref{J32}.  Therefore, its homotopy class
in~$\Bolza^{(1)}$ can be represented by an imbedded loop $\beta\subset
\Bolza^{(1)}$ of length at most~$8y$.  Thus,~$\beta$ contains fewer
than 8 edges of $\Bolza^{(1)}$.  Since the number of edges must be
even, its image under~$Q$ must retract to a circuit with at most six
edges in the 1-skeleton of the cubical subdivision of~$S^2$.  If the
number is 4, then the circuit lies in the boundary of a square face of
the cube in $S^2$.  But the boundary of a face does not lift
to~$\Bolza$, since it surrounds a single ramification point, namely
the center of the square face.

Hence there must be 6 edges in the circuit.  There are two types of
circuits with six edges in the 1-skeleton of the cubical subdivision
of~$S^2$:
\begin{enumerate}
\item[(a)]
the boundary of the union of a pair of adjacent squares; 
\item[(b)] a path consisting of the edges meeting a suitable great
circle.
\end{enumerate}
However, a path of type (b) surrounds an odd number, namely 3, of
ramification points, and hence does not lift to the genus 2 surface.
Meanwhile, a path of type (a) surrounds two ramification points, and
hence does lift to the surface.  Such a path is freely homotopic
in~$\Bolza$ to one of the 12 geodesics of type
$\gamma$ (\cf Lemma~\ref{J21}), completing the proof.
\end{proof}

\section{Voronoi cells and Euler characteristic}
\label{Jproof12}

The following proposition provides a preliminary lower bound on the
area of hyperelliptic surfaces with nonpositive curvature.

\begin{proposition} 
\label{Jprop:nonpositive}
Every~$J$-invariant~$\CAT(0)$ metric~$\gmetric$ on a closed
hyperelliptic surface~$\Sigma_\genus$ of genus~$\genus$, satisfies the
bound 
\[
\SR(\Sigma_\genus, \gmetric) \leq {8} ({(\genus+1)\pi})^{-1} .
\]
\end{proposition}

\begin{proof}
To prove this scale-invariant inequality, we normalize the metric
on~$\Sigma=\Sigma_\genus$ to unit systole,
\ie~$\pisys_1(\Sigma,\gmetric)=1$.  The preimage by~$Q:\Sigma \to S^2$
of an arc of~$S^2$ joining two distinct branch points forms a
noncontractible loop on~$\Sigma$.  Therefore, the distance between two
Weierstrass points is at least~$\frac{1}{2}
\pisys_1(\Sigma,\gmetric) = \frac{1}{2}$. Thus, we obtain $2\genus+2$
disjoint disks of radius~$\frac{1}{4}$, centered at the Weierstrass
points.  Since the metric is~$\CAT(0)$, the area of each disk is at
least~$\frac{\pi}{16}$.  Thus,~$\area(\Sigma,\gmetric) \geq
\frac{\genus+1}{8} \pi$.
\end{proof}

An {\em optimal\/} lower bound requires a more precise estimate on the
area of the Voronoi cells.  The idea is to replace area of balls by
area of polygons, where control over the number of sides is provided
by the Euler characteristic, \cf \cite{Bav2}.

Denote by~$u: \tilde{\Sigma} \rightarrow \Sigma$ its universal cover.  Let
$\{x_i \mid i \in \N \}$ be an enumeration of the lifts of Weierstrass
points on~$\tilde{\Sigma}$.  The Voronoi cell~$V_i \subset \tilde \Sigma$
centered at~$x_i$ is defined as the set of points closer to~$x_i$ than
to any other lift of a Weierstrass point.  In formulas,
\begin{equation}
\label{J}
V_i = \left\{ \left. x \in \tilde{\Sigma} \right| d(x,x_i) \leq d(x,x_j)
\text{ for every } j \neq i \right\}.
\end{equation}
The Voronoi cells on~$\tilde{\Sigma}$ are polygons whose edges are arcs of
the equidistant curves between a pair of lifts of Weierstrass
points. Note that these edges are not
necessarily geodesics.  The Voronoi cells on~$\tilde{\Sigma}$ are
topological disks, while their projections~$u(V_i)\subset \Sigma$ may have
more complicated topology.  Thus, the surface~$\Sigma$ decomposes
into~$2\genus+2$ images of Voronoi cells, centered at the~$2\genus+2$
Weierstrass points.  By the number of sides of~$u(V)$ we will mean the
number of sides of the polygon~${V}$.

\begin{lemma}
Let $\gmetric_{1}$ and $\gmetric_{2}$ be two $\CAT(0)$ metrics lying
in the same conformal class.  Then, the averaged metric
$\gmetric=\tfrac{1}{2} (\gmetric_{1}+\gmetric_{2})$ is $\CAT(0)$, as
well.
\end{lemma}

\begin{proof}
Choose a point~$x \in \Sigma$ and a metric~$\gmetric_{0}$ in its
conformal class, which is flat in a neighborhood of~$x$.  Every metric
$\gmetric = H \gmetric_{0}$ conformal to it satisfies
\[
K_{\gmetric} H = K_{\gmetric_{0}} - \tfrac{1}{2} \Delta \log H,
\]
\cf \cite[p.~265]{GHL}, where $K_{\gmetric}$ and $K_{\gmetric_{0}}$
are the Gaussian curvatures of $\gmetric$ and~$\gmetric_{0}$, and
$\Delta$ is the Laplacian of~$\gmetric_{0}$ with $\Delta f =
\mbox{div} \nabla f$.  Thus, the metrics $\gmetric_{i}$ can be written
as $\gmetric_{i} = e^{h_{i}} \gmetric_{0}$ where $h_{i}$ is
subharmonic in the neighborhood of~$x$, \ie $\Delta h_{i} \geq 0$.  A
simple computation shows that
\[
\Delta \log H \geq \frac{e^{h_{1}} \Delta h_{1} + e^{h_{2}} \Delta
h_{2}}{2H} \geq 0
\]
where $H = \tfrac{1}{2} (e^{h_{1}} + e^{h_{2}})$, proving the lemma if
both points are regular.  For singular points with positive angle
excess, the CAT(0) property for the averaged metric is immediate from
Remark~\ref{tro}.
\end{proof}

\begin{proof}[Proof of Theorem~\ref{Jtheo:nonpositive}]
Since averaging by $J$ can only improve the systolic ratio, we may
assume without loss of generality that our metric is already
$J$-invariant.

\forget
Consider the exponential map~$\exp_A$ at a Weierstrass point~$A\in \Sigma$.
If~$A$ is a singular point, we lift to the appropriate smooth finite
cover, and use its exponential map, \cf Definition~\ref{911}.

Let~$B\in \Sigma$ be another Weierstrass point.  Consider its lift~$B_0 =
\exp_A^{-1}(B)$ to the tangent plane~$T_A$.  Consider the equidistant
line
\[ 
L_{O, B_0} \subset T_A
\]
between the origin~$O\in T_A$ and the point~$B_0$.

By the Toponogov theorem applied to the hinge at~$O$ defined by the
point~$B_0$ and an arbitrary point~$X\in L_{O, B_0}$ on the
equidistant line, we have
\begin{equation}
\label{J44b}
\dist_\Sigma(A, \exp_A(X)) \leq \dist_\Sigma(B, \exp_A(X)) .
\end{equation}
\forgotten

There exists an extension of the notions of tangent plane and
exponential map, to surfaces with singularities.  Namely, let $A\in
\Sigma$.  There exists a CAT(0) piecewise flat plane~$T_{A}$ with
conical singularities and a covering $\exp_{A}:T_{A} \to \Sigma$ with
the following properties:
\begin{enumerate}
\item 
$\exp_A$ sends the origin $O$ of $T_{A}$ to~$A$;

\item 
$\exp_A$ takes the conical singularities of~$T_{A}$ to the
singularities of~$\Sigma$;

\item 
$\exp_A$ sends every pair of geodesic arcs issuing from the origin $O
\in T_{A}$ and forming an angle~$\theta$, to a pair of geodesic arcs
of the same lengths, and forming the same angle, at their
basepoint~$A$.
\end{enumerate}

By the Rauch Comparison Theorem, the exponential map $\exp_{A}$ does
not decrease distances.

Now assume $A\in \Sigma$ is a Weierstrass point, and let~$B\in \Sigma$
be another Weierstrass point.  Fix a lift~$B_0$ of~$B$ to the tangent
plane~$T_A$, along a minimizing arc.  Consider the equidistant line
\[ 
L_{O, B_0} \subset T_A
\]
between the origin~$O\in T_A$ and the point~$B_0$.  Consider a
point~$X_{0}\in L_{O, B_0}$.  Let $X=\exp_{A}(X_{0})$.  Since the
exponential map does not decrease distances, we have
\begin{equation}
\label{J44b}
\dist_\Sigma(A, X) = \dist_{T_{A}}(O,X_{0}) = \dist_{T_{A}} (B_{0},
X_{0}) \leq \dist_\Sigma(B, X) .
\end{equation}

Now consider the polygon in the tangent plane~$T_{A}$, obtained as the
intersection of the halfspaces containing the origin, defined by the
lines~$L_{O, B_0}$, as~$B$ runs over all Weierstrass points.  It
follows from~\eqref{J44b} that the exponential image of this polygon
is contained in the Voronoi cell of~$A$.  Since the exponential map
does not decrease distances, the area of the polygon is a lower bound
for the area of the Voronoi cell.  If~$k$ is the number of sides
of~$V$, then~$V$ is partitioned into~$k$ triangles with
angle~$\theta_i$ at~$O$, whose area is bounded below by
\begin{equation}
\label{J43}
\left( \tfrac{\pisys_1}{4} \right)^2 \tan\tfrac{\theta_i}{2}
\end{equation}
since $\dist(A,B) \geq \frac{\pisys_{1}}{2}$.

Consider the graph on~$S^2$ defined by the projections of the Voronoi
cells to the sphere.  Thus we have~$f=6$ faces.  Applying the formula
$v-e+f=2$ and the well-known fact that~$3v \leq {2} e$, we obtain
\[
{e} \leq 3f-6 =12.
\]
Hence the spherical graph has at most~12 edges.  Note that the maximum
is attained by the 1-skeleton of the cubical subdivision.

The area of a flat isosceles triangle with third angle~$\theta$, and
with unit altitude from the third vertex, is~$\tan \tfrac {\theta}
{2}$.  This formula provides a lower bound for the area of the Voronoi
cells as in \eqref{J43}.  The proof is completed by Jensen's
inequality applied to the convex function~$\tan (\tfrac {x} {2})$
when~$0 < x < \pi$.  In the boundary case of equality, we have~$e=12$,
all angles~$\theta_i$ as in~\eqref{J43} must be equal, curvature must
be zero because of equality in the Rauch Comparison Theorem, and we easily deduce that
each Voronoi cell is a regular octagon.  To minimize the area of the
octagon, we must choose~$\theta$ as small as possible.  The~$\CAT(0)$
hypothesis at the center of the octagon imposes a lower bound~$\theta
\geq \tfrac \pi{4}$.  Hence the optimal systolic ratio is achieved for
the regular flat octagon with a smooth point at the center.
\end{proof}

\section{Arbitrary metrics on the Bolza surface}
\label{Jfive}

The conformal class of the Bolza surface~$\Bolza$, \cf
Section~\ref{J22}, is likely to contain a systolically optimal surface
in genus~2, as discussed in Section~\ref{Jone1}.

\begin{theorem} \label{theo:conf}
Every metric~$\gmetric$ in the conformal class of the Bolza surface
satisfies the bound
\begin{equation}
\label{J12}
\SR(\gmetric) \leq \tfrac{\pi}{3} = 1.0471\ldots
\end{equation}
\end{theorem}

\begin{remark}
In particular, every metric~$\gmetric$ in the conformal class of the
Bolza surface satisfies Bavard's inequality
\begin{equation}
\label{Jcb}
\SR(\gmetric) \leq \frac{\pi}{2^{3/2}} \simeq 1.1107
\end{equation}
for the Klein bottle \cite{Bav1}.  This suggests a possible
monotonicity of~$\chi(\Sigma)$ as a function of~$\SR(\Sigma)$.
\end{remark}

\begin{lemma}
\label{J54}
Let~$\gmetric$ be an~$\Aut(\Bolza)$-invariant metric on~$\Bolza$.
Consider the displacement~$\delta(J)$ of~$J$ on the
1-skeleton~$\Bolza^{(1)}$ of the Voronoi subdivision of~$\Bolza$,
namely
\[
\delta(J) = \min_{x\in \Bolza^{(1)}} \dist(x, J(x)).
\]
Then~$\area(\Bolza,\gmetric) \geq 6 \left( \tfrac{2}{\pi} \right)
\delta(J)^2$.
\end{lemma}

\begin{proof}
Consider the Voronoi subdivision with respect to the set of six
Weierstrass points on~$\Bolza$.  Since each Voronoi
cell~$\Omega\subset \Bolza$ is~$J$-invariant, we can identify all
pairs of opposite points of the boundary~$\partial \Omega$, to obtain
a projective plane
\begin{equation}
\label{J72} 
\rp^2 = \Omega/\sim.
\end{equation}
We now apply Pu's inequality \cite{Pu} to each of the six Voronoi
cells, to obtain
\begin{equation}
\label{J62}
\area(\rp ^2)\geq \tfrac{2}{\pi} \pisys_1(\rp^2)^2 .
\end{equation}
The lemma now follows from the bound~$\pisys_1(\rp^2) \geq \delta(J)$.
\end{proof}

\begin{lemma}
\label{J74}
Every~$\Aut(\Bolza)$-invariant metric~$\gmetric$ on~$\Bolza$ satisfies
the bound $2 \delta(J) \geq \pisys_1(\Bolza, \gmetric)$.
\end{lemma}

\begin{proof}
Let~$p, J(p)\in \partial \Omega$ satisfy~$\dist(p, J(p)) = \delta(J)$.
Let~$\alpha \subset \Omega$ be a minimizing path joining~$p$
to~$J(p)$.  Let~$\Omega'\subset \Bolza$ be the adjacent Voronoi cell
containing this pair of boundary points, and~$r: \Bolza \to \Bolza$
the anticonformal involution which switches~$\Omega$ and~$\Omega'$,
and fixes their common boundary.  The loop~$\alpha \cup r(\alpha)$
belongs to the free homotopy class of the noncontractible loop~$\gamma
\subset \Bolza$ obtained as the inverse image under~$Q: \Bolza \to
S^2$ of the edge of the octahedral decomposition of~$S^{2}$ joining
the images of the centers of~$\Omega$ and~$\Omega'$, \cf
Lemma~\ref{J21}.  Since~$\length(\alpha \cup r(\alpha)) = 2
\delta(J)$, the lemma follows.
\end{proof}

\begin{proof}[Proof of Theorem~\ref{theo:conf}]
We may assume that the metric on~$\Bolza$ is~$\Aut(\Bolza)$-invariant,
since averaging the metric by a finite group of holomorphic and
antiholomorphic diffeomorphisms can only improve the systolic ratio.
We combine the inequalities of Lemmas~\ref{J74} and \ref{J54} to prove
the theorem.
\end{proof}



\begin{thebibliography}{ABCDe}

\bibitem[Am04]{Am} Ammann, B.: Dirac eigenvalue estimates on two-tori.
{\em J. Geom. Phys.} \textbf{51} (2004), no. 3, 372--386.


\bibitem[BCIK05]{BCIK1} Bangert, V; Croke, C.; Ivanov, S.; Katz, M.:
Filling area conjecture and ovalless real hyperelliptic surfaces, {\em
Geometric and Functional Analysis (GAFA)\/} \textbf{25} (2005), no.~3.
577-597.  See \texttt{arXiv:math.DG/0405583}

\bibitem[BCIK06]{BCIK2} Bangert, V; Croke, C.; Ivanov, S.; Katz, M.:
Boundary case of equality in optimal Loewner-type inequalities, {\em
Trans. Amer. Math. Soc.} \textbf{358} (2006).  See
\texttt{arXiv:math.DG/0406008}



\bibitem[BK03]{BK1} Bangert, V.; Katz, M.: Stable systolic
inequalities and cohomology products, {\em Comm. Pure Appl. Math.}
\textbf{56} (2003), 979--997.  Available at
\texttt{arXiv:math.DG/0204181}

\bibitem[BK04]{BK2} Bangert, V; Katz, M.: An optimal Loewner-type
systolic inequality and harmonic one-forms of constant norm.  {\em
Comm. Anal. Geom.}  \textbf{12} (2004), number 3, 703-732.  See
\texttt{arXiv:math.DG/0304494}

\bibitem[Bav86]{Bav1} Bavard, C.: In\'egalit\'e isosystolique pour la
bouteille de Klein.  {\em Math. Ann.} \textbf{274} (1986), no. 3,
439--441.


\bibitem[Bav92]{Bav2} Bavard, C.: La systole des surfaces
hyperelliptiques, {\em Prepubl. Ec. Norm. Sup. Lyon\/} \textbf{71}
(1992).



\bibitem[BH99]{BH} Bridson, M.; Haefliger, A.: Metric spaces of
non-positive curvature.  {\em Grundlehren der Mathematischen
Wissenschaften}, \textbf{319}.  Springer-Verlag, Berlin, 1999.




\bibitem[Br96]{Bry} Bryant, R.: On extremals with prescribed
Lagrangian densities, Manifolds and geometry (Pisa, 1993), 86--111,
{\em Sympos. Math.}, \textbf{36}, Cambridge Univ. Press, Cambridge,
1996.

\bibitem[BS94]{BS} Buser, P.; Sarnak, P.: On the period matrix of a
Riemann surface of large genus.  With an appendix by J. H. Conway and
N. J. A. Sloane.  {\em Invent. Math.} \textbf{117} (1994), no.~1,
27--56.

\bibitem[Ca96]{Ca} Calabi, E.: Extremal isosystolic metrics for
compact surfaces, Actes de la table ronde de geometrie differentielle,
{\em Sem. Congr.}  \textbf{1} (1996), Soc. Math. France, 167--204.


\bibitem[CK03]{CK} Croke, C.; Katz, M.: Universal volume bounds in
Riemannian mani\-folds, {\it Surveys in Differential Geometry\/}
\textbf{8}, Lectures on Geometry and Topology held in honor of Calabi,
Lawson, Siu, and Uhlenbeck at Harvard University, May 3- 5, 2002,
edited by S.T. Yau (Somerville, MA: International Press, 2003.)
pp. 109 - 137.  See \texttt{arXiv:math.DG/0302248}


\bibitem[FK92]{FK} Farkas, H.; Kra, I.: Riemann surfaces.  Second
edition. Graduate Texts in Mathematics, 71. Springer-Verlag, New York,
1992.

\bibitem[GHL90]{GHL} Gallot, S.; Hulin, D.; Lafontaine, J.: Riemannian
Geometry. Springer-Verlag, Berlin, 1990.

\bibitem[Gr83]{Gr1} Gromov, M.: Filling Riemannian manifolds, {\em
J. Diff. Geom.\/} \textbf{18} (1983), 1-147.

\bibitem[Gr96]{Gr2} Gromov, M.: Systoles and intersystolic
inequalities.  Actes de la Table Ronde de G\'{e}om\'{e}trie
Diff\'{e}rentielle (Luminy, 1992), 291--362, {\em S\'{e}min. Congr.}
\textbf{1}, Soc. Math. France, Paris, 1996.
\newline\noindent
\texttt{www.emis.de/journals/SC/1996/1/ps/smf\_sem-cong\_1\_291-362.ps.gz}

\bibitem[Gr99]{Gr3} Gromov, M.: Metric structures for Riemannian and
non-Riemannian spaces.  {\em Progr. in Mathematics}, \textbf{152},
Birkh\"{a}user, Boston, 1999.

\bibitem[IK04]{IK} Ivanov, S.; Katz, M.: Generalized degree and
optimal Loewner-type inequalities, {\em Israel J. Math.}  \textbf{141}
(2004), 221-233.  \texttt{arXiv:math.DG/0405019}


\bibitem[JP04]{JP} Jakobson, D.; Levitin, M.; Nadirashvili, N.;
Polterovich, I.: Spectral problems with mixed Dirichlet-Neumann
boundary conditions: isospectrality and beyond, preprint.  See
\texttt{arXiv:math.SP/0409154}

\bibitem[Je84]{Je} Jenni, F.: \"Uber den ersten Eigenwert des
Laplace-Operators auf ausgew\"{a}hlten Beispielen kompakter
Riemannscher Fl\"{a}chen. [On the first eigenvalue of the Laplace
operator on selected examples of compact Riemann surfaces] {\em
Comment. Math. Helv.} \textbf{59} (1984), no. 2, 193--203.


\bibitem[Ka03]{Ka3} Katz, M.: Four-manifold systoles and surjectivity
of period map, {\em Comment. Math. Helv.} \textbf{78} (2003), 772-786.
\texttt{arXiv:math.DG/0302306}

\bibitem[Ka06]{Ka4} Katz, M.: Systolic inequalities and Massey
products in simply-connected manifolds.  See
\texttt{arXiv:math.DG/0604012}


\bibitem[KL05]{KL} Katz, M.; Lescop, C.: Filling area conjecture,
optimal systolic inequalities, and the fiber class in abelian covers.
Geometry, spectral theory, groups, and dynamics, 181--200, {\em
Contemp. Math.} \textbf{387}, Amer. Math. Soc., Providence, RI, 2005.
See \texttt{arXiv:math.DG/0412011}


\bibitem[KR06]{KR} Katz, M.; Rudyak, Y.: Lusternik-Schnirelmann
category and systolic category of low dimensional manifolds.  {\em
Communications on Pure and Applied Mathematics}, \textbf{59} (2006).
See \texttt{arXiv:math.DG/0410456}

\bibitem[KR07]{KR2} Katz, M.; Rudyak, Y.: Bounding volume by systoles
of 3-manifolds.  See \texttt{arXiv:math.DG/0504008}


\bibitem[KRS07]{KRS} Katz, M.; Rudyak, Y.: Sabourau, S.: Systoles of
2-complexes, Reeb graph, and Grushko decomposition.  See
\texttt{arXiv:math.DG/0602009}


\bibitem[KS05]{KS2} Katz, M.; Sabourau, S.: Entropy of systolically
extremal surfaces and asymptotic bounds, {\em Ergodic Theory and
Dynamical Systems\/} \textbf{25} (2005), no.~4, 1209-1220.  Available
at \texttt{arXiv:math.DG/0410312}

\bibitem[KS06]{KS1} Katz, M.; Sabourau, S.: Hyperelliptic surfaces are
Loewner, {\em Proc. Amer. Math. Soc.} \textbf{134} (2006), no.~4,
1189-1195.  See \texttt{arXiv:math.DG/0407009}

\bibitem[KSV06]{KSV} Katz, M.; Schaps, M.; Vishne, U.: Logarithmic
growth of systole of arithmetic Riemann surfaces along congruence
subgroups.  Available at \texttt{arXiv:math.DG/0505007}

\bibitem[Kon03]{Kon} Kong, J.: Seshadri constants on Jacobian of
curves.  {\em Trans. Amer. Math. Soc.} \textbf{355} (2003), no. 8,
3175--3180.

\bibitem[KuN95]{KuN} Kuusalo, T.; N\"a\"at\"anen, M.: Geometric
uniformization in genus~$2$.  {\em Ann. Acad. Sci. Fenn. Ser. A Math.}
\textbf{20} (1995), no. 2, 401--418.

\bibitem[Mi95]{Mi} Miranda, R.: Algebraic curves and Riemann
surfaces. {\em Graduate Studies in Mathematics}, {\textbf 5}.
Amer. Math. Soc., Providence, RI, 1995.

\bibitem[Pu52]{Pu} Pu, P.M.: Some inequalities in certain
nonorientable Riemannian manifolds, {\it Pacific J. Math.\/}
\textbf{2} (1952), 55--71.


\bibitem[RS07]{RS} Rudyak, Y.: Sabourau, S.: Systolic invariants of
groups and $2$-complexes via Grushko decomposition.

\bibitem[Sa04]{Sa1} Sabourau, S.: Systoles des surfaces plates
singuli\`eres de genre deux, {\em Math. Zeitschrift\/} \textbf{247}
(2004), no. 4, 693--709.

\bibitem[Sab06]{Sa2} Sabourau, S.: Systolic volume and minimal entropy
of aspherical manifolds, {\em J. Differential Geom.}, to appear.  See
\texttt{arXiv:math.DG/0603695}

\bibitem[Sa06b]{Sa3} Sabourau, S.: Asymptotic bounds for separating
systoles on surfaces, preprint.

\bibitem[Sc93]{Sc1} Schmutz, P.: Riemann surfaces with shortest
geodesic of maximal length.  {\em Geom. Funct. Anal.} \textbf{3}
(1993), no. 6, 564--631.

\bibitem[Tr90]{Tro} Troyanov, M.: Coordonn\'ees polaires sur les
surfaces riemanniennes singuli\`eres.  {\em Ann. Inst. Fourier
(Grenoble)} \textbf{40} (1990), no. 4, 913--937 (1991).

\end{thebibliography}
\end{document}